\documentclass[twoside,11pt]{article} 
\usepackage{amsfonts} 
\usepackage{graphicx}
\usepackage{color}

\newtheorem{thm}{Theorem}

\begin{document}

 \def\reals{{\mathbb R}}
 \def\ch{{\cal H}}
 \def\cA{{\cal A}}
 \def\cD{{\cal D}}
 \def\cK{{\cal K}}
 \def\cC{{\cal C}}
 \def\cN{{\cal N}}
 \def\cR{{\cal R}}
 \def\cS{{\cal S}}
 \def\cT{{\cal T}}
 \def\cV{{\cal V}}
 %notation for cabling torus
 \def\ta{{\cal T}_{\subset}}
 \def\cI{{\cal I}}
 %notation for cabling operation
 \def\bC{{\bf C}}
 %notation for braid axis
 \def\axis{{\bf A}}
 %notation for braid fibration
 \def\fibr{{\bf H}}
 %a-arcs
 \def\ba{{\bf a}}
 %b-arcs
 \def\bb{{\bf b}}
 %c-arcs
 \def\bc{{\bf c}}
 %e-arcs
 \def\be{{\bf e}}
 \def\d{{\delta}} 
 \def\ci{{\circ}} 
 \def\e{{\epsilon}} 
 \def\l{{\lambda}} 
 \def\L{{\Lambda}} 
 \def\m{{\mu}} 
 \def\n{{\nu}} 
 \def\o{{\omega}} 
 \def\s{{\sigma}} 
 \def\v{{\varphi}} 
 \def\a{{\alpha}} 
 \def\b{{\beta}} 
 \def\p{{\partial}} 
 \def\r{{\rho}} 
 \def\ra{{\rightarrow}} 
 \def\lra{{\longrightarrow}} 
 \def\g{{\gamma}} 
 \def\D{{\Delta}} 
 \def\La{{\Leftarrow}} 
 \def\Ra{{\Rightarrow}} 
 \def\x{{\xi}} 
 \def\c{{\mathbb C}} 
 \def\z{{\mathbb Z}} 
 \def\2{{\mathbb Z_2}} 
 \def\q{{\mathbb Q}} 
 \def\t{{\tau}} 
 \def\u{{\upsilon}} 
 \def\th{{\theta}} 
 \def\la{{\leftarrow}} 
 \def\lla{{\longleftarrow}} 
 \def\da{{\downarrow}} 
 \def\ua{{\uparrow}} 
 \def\nwa{{\nwtarrow}} 
 \def\swa{{\swarrow}} 
 \def\nea{{\netarrow}} 
 \def\sea{{\searrow}} 
 \def\hla{{\hookleftarrow}} 
 \def\hra{{\hookrightarrow}} 
 \def\sl{{SL(2,\mathbb C)}} 
 \def\ps{{PSL(2,\mathbb C)}} 
 \def\qed{{\hfill$\diamondsuit$}} 
 \def\pf{{\noindent{\bf Proof.\hspace{2mm}}}} 
 \def\ni{{\noindent}} 
 \def\sm{{{\mbox{\tiny M}}}} 
 \def\sc{{{\mbox{\tiny C}}}} 
 \def\ke{{\mbox{ker}(H_1(\p M;\2)\ra H_1(M;\2))}} 
 \def\et{{\mbox{\hspace{1.5mm}}}}

\markboth{William W. Menasco \today}{Iterated Torus Knots \today}
\title{Erratum:  On iterated torus knots and transversal knots} 
 
\author{William W. Menasco \thanks{partially supported by NSF grant \#DMS 0306062}}  

\date{\today}

\maketitle

In \cite{[M1]} the author's Theorems 1.1 and 1.2,  combined, implied that iterated torus knots are transversally simple.
The previous best result in this direction had been the proof, due to John Etnyre, that positive torus knots are transversally simple.  

After the publication of \cite{[M1]}, Etnyre and Honda argued in
\cite{[EH]}, using very different techniques from those in \cite{[M1]}, that among the infinitely many Legendrian knots whose topological knot type is that of the $(3,2)$-cable on the
$(3,2)$-torus knot,  there must be one whose transversal pushoff cannot be transversally simple.  This lead to a review of the proof in \cite{[M1]} and the discovery of an error.  

To pinpoint the error, we recall that in \cite{[BW]} it was established that exchange reducible
closed braids representing oriented knots are transversally simple when the braids are
viewed as transversal knots in the standard contact structure of $\reals^3$ or $S^3$.
The author's Theorem 1.1 of \cite{[M1]} asserted that all iterated torus knots are exchange reducible. Theorem 1.2  then followed.  The error occurs in the proof of Lemma 5.13 of \cite{[M1]}, needed for the proof of Theorem 1.1.   The changes made in the passage from Figure 15b to 16b is not valid, due to an obstruction which will be described in detail in a manuscript now in preparation.  In said manuscript we will be able  to recover the major part of the claimed results, and to show exactly why and how the examples pointed to in \cite{[EH]} and other related examples fail to be transversally simple.  

Without introducing major new techniques, we can only claim a simplified version of Theorems 1.1 and 1.2 of \cite{[M1]}:

\begin{thm}
\label{theorem:torus knots are simple}
Torus knots are exchange reducible, and hence transversally simple.
\end{thm}

A deeper understanding of the error in Lemma 5.13 yields a complete understanding of which iterated torus
knots are exchange reducible and is investigated in \cite{[M2]}

\end{document}